\documentclass[10pt,twoside]{article}

\usepackage{graphicx}
\usepackage{mathtools}
\usepackage [usenames] {color}
\usepackage [colorlinks=true,
linkcolor=webgreen,
filecolor=webbrown,
citecolor=webgreen] {hyperref}
\definecolor {webgreen} {rgb} {0,.5,0}

\definecolor {webbrown} {rgb} {.6,0,0}
      
\usepackage {color}
\usepackage {float}

\usepackage{amsfonts}
\usepackage{amstext}
\usepackage{amssymb,amsmath,amsbsy}

\usepackage{hyperref}

\newtheorem{theorem}{Theorem}
\newtheorem{problem}{Problem}

\newtheorem{definition}{Definition}
\newtheorem{lemma}{Lemma}
\newtheorem{proposition}{Proposition}

\newcounter{thcount}
\def\0{\mathbf{0}}
\def\1{\mathbf{1}}

\def\X{\mathbf{X}}

\def\C{\mathbf{C}}

\def\A{\mathbf{A}}

\def\DD{\mathbf{D}}

\def\diag{\mathrm{diag}}

\def\rk{\mathrm{rank}}

\def\Vol{{\mathrm{Vol}}}

\def\disc{{\mathrm{disc}}}
\def\part{\cal P} 

\begin{document}


{\Large {\bf {Properties of the multiway discrepancy}}\\

{\large {\bf Marianna Bolla$^{1,2}$, Edward Kim$^{2,3}$, Cheng Wai Koo$^{2,4}$}}\\

\begin{center}
$^1$Institute of Mathematics, Budapest University of Technology and Economics,
Hungary
\end{center}

\begin{center}
$^2$Budapest Semester of Mathematics, Hungary
\end{center}

\begin{center}
$^3$Duke University, NC, USA
\end{center}

\begin{center}
$^4$Harvey Mudd College, CA, USA
\end{center}


\vskip1cm

\renewcommand\abstractname{Abstract}

\begin{abstract}
\noindent \\
Some properties of the multiway discrepancy of rectangular matrices of
nonnegative entries are discussed. We are able to prove the continuity of 
this discrepancy, as well as some statements about the multiway discrepancy
of some special matrices and graphs. We also conjecture that the $k$-way
discrepancy  is monotonic in $k$.

\noindent
\textbf{Keywords:} {Multiway discrepancy; Edge-density; Binary arrays; Independent table.}

\end{abstract}

\section {Introduction}
\label{intro}

In many applications, for example when microarrays are analyzed, our
data are collected in the form of an $m\times n$ rectangular array 
$\A=(a_{ij})$ of
nonnegative real entries, called contingency table. 
We assume that $\A$ is
non-decomposable, i.e., $\A \A^T$ (when $m\le n$) or 
$\A^T \A$  (when $m > n$) is irreducible.
Consequently,
the row-sums
$d_{row,i} =\sum_{j=1}^n a_{ij}$ and column-sums $d_{col,j}=\sum_{i=1}^m a_{ij}$
of $\A$ are strictly positive, and the diagonal matrices
$\DD_{row} =\diag (d_{row,1} ,\dots ,d_{row,m})$ and 
$\DD_{col} =\diag (d_{col,1} ,\dots ,d_{col,n})$ are regular.
Without loss of generality, we also assume that
$\sum_{i=1}^n \sum_{j=1}^m a_{ij} =1$, since the
normalized table 
\begin{equation}\label{cnor}
 \A_{D} =  \DD_{row}^{-1/2} \A  \DD_{row}^{-1/2} ,
\end{equation}
is not affected by the scaling 
of the entries of $\A$.
It is well known (see e.g.,~\cite{Bolla14}) that the singular values  of 
$\A_{D}$ are in the [0,1]
interval. Enumerated in non-increasing order, they are the real numbers
$$
 1=s_0 >s_1 \ge \dots \ge s_{r-1} > s_{r} = \dots = s_{n-1} =0  ,
$$
where $r= \rk (\A )$. When $\A$ is non-decomposable,
1 is a single singular value, and it is
denoted by $s_0$, since it belongs to the trivial singular vector pair.
In~\cite{Bollarx} we gave an upper estimate for $s_k$ in terms of the
$k$-way discrepancy to be introduced herein. In~\cite{Bolla14} a certain
converse of this estimate was proved. Therefore, by the monotonic property
of the singular values we guess, that the $k$-way discrepancy is also monotonic
decreasing in $k$.

\begin{definition}\label{diszkrepancia}
The multiway discrepancy  of the rectangular array $\A$ of nonnegative entries
in the proper $k$-partition $R_1 ,\dots ,R_k$ of its rows and
$C_1 ,\dots ,C_k$ of its columns is
\begin{equation}\label{disk} 
 \disc (\A ; R_1 ,\dots ,R_k , C_1 ,\dots ,C_k ) =
 \max_{\substack{1\le a , b\le k \\X\subset R_a , \, Y\subset C_b}} 
 \disc (X,Y;R_a , C_b ) ,
\end{equation}
where for $X\subset R_a $ and $Y\subset C_b$,
\begin{equation}
\begin{aligned}
\disc (X,Y;R_a , C_b ) &=
 \frac{|a (X, Y)-\rho (R_a,C_b ) \Vol (X)\Vol (Y)|}{\sqrt{\Vol(X)\Vol(Y)}} \\
 &=| \rho (X,Y) -\rho (R_a,C_b ) | \sqrt{\Vol(X)\Vol(Y)}.
\end{aligned}
\end{equation}
Here
$a (X, Y) =\sum_{i\in X} \sum _{j\in Y} a_{ij}$ is the cut between
$X\subset R_a$ and $Y\subset C_b$, 
$\Vol (X) = \sum_{i\in X} d_{row,i}$ is the volume of the row-subset $X$, 
$\Vol (Y) = \sum_{j\in Y} d_{col,j}$ is the volume of the column-subset $Y$, 
whereas
$\rho (X,Y) =\frac{a(X,Y)}{ \Vol (X) \Vol (Y)}$ denotes the 
density between $X$ and $Y$.
The minimum $k$-way discrepancy  of  $\A$ itself is
$$
 \disc_k (\A ) = \min_{\substack{R_1 ,\dots ,R_k \\ C_1 ,\dots ,C_k } } 
 \disc (\A ; R_1 ,\dots ,R_k , C_1 ,\dots ,C_k ).
$$
\end{definition}

\section{Continuity of the discrepancy}
\label{cont}

\begin{theorem}[Continuity of the density]
Let $\A$ be an $m \times n$ binary array where each entry is $0$ or $1$, with row set $R$ and column set $C$. Define the \emph{volume} of any subset $X \subset R$ or any subset $Y \subset C$ as
\[
V(X) = \sum_{i \in X} \sum_{j \in C} a_{ij}, \qquad V(Y) = \sum_{i \in R} \sum_{j \in Y} a_{ij}.
\]
Assume no row and no column has only zeros, so $V(X), V(Y) \ge 1$ for any nonempty $X \subset R$ or $Y \subset C$. For any nonempty subsets $X \subset R$, $Y \subset C$, define
\[
d(X,Y) = \frac{a(X,Y)}{V(X) \cdot V(Y)}\,,
\]
where
\[
a(X,Y) = \sum_{i\in X}\sum_{j\in Y}a_{ij}.
\]
Then, for any $X^{*}\subset X$ with $\frac{V(X^{*})}{V(X)} \ge 1-\delta$ and any $Y^{*}\subset Y$ with $\frac{V(Y^{*})}{V(Y)} \ge 1-\delta$, this Theorem states that
\[
|d(X,Y) - d(X^*, Y^*)| \le 4\delta
\]
for $0< \delta <\frac{1}{4}$.
\end{theorem}

\noindent
\textbf{Proof.}
Since $a(X^*,Y^*) \le a(X,Y)$,
\[
d(X^*,Y^*) = \frac{a(X^*,Y^*)}{V(X^*) \cdot V(Y^*)} \le \frac{a(X,Y)}{V(X) \cdot V(Y)}\cdot \frac{V(X)}{V(X^*)} \cdot \frac{V(Y)}{V(Y^*)}\,.
\]
Since $\frac{a(X,Y)}{V(X) \cdot V(Y)} = d(X,Y)$, $\frac{V(X)}{V(X^*)} \le \frac{1}{1-\delta}$, and $\frac{V(Y)}{V(Y^*)} \le \frac{1}{1-\delta}$,
\[
d(X^*,Y^*) \le d(X,Y)\cdot \left(\frac{1}{1-\delta}\right)^2, \tag{1}
\]
which gives an upper bound for $d(X^*,Y^*)$. For the lower bound, observe that
\[
a(X^*,Y^*) = a(X,Y) - [a(X^*,Y-Y^*)+a(X-X^*,Y^*)+a(X-X^*,Y-Y^*)] \tag{2}
\]
Since $0 \le d(X,Y) \le 1$ (it is true only in this binary case!), for any $X \subset R$, $Y \subset C$,
\begin{align*}
a(X^*, Y-Y^*) &\le V(X^*) \cdot \left( V(Y) - V(Y^*) \right), \tag{3.1}\\
a(X-X^*, Y^*) &\le \left( V(X) - V(X^*) \right) \cdot V(Y^*), \tag{3.2}\\
a(X-X^*, Y-Y^*) &\le \left( V(X) - V(X^*) \right) \cdot \left( V(Y) - V(Y^*) \right). \tag{3.3}
\end{align*}\\
Dividing both sides of (2) by $V(X^*)\cdot V(Y^*)$ and applying (3.1)-(3.3),
\begin{align*}
d(X^*,Y^*) &\ge \frac{a(X,Y)}{V(X)\cdot V(Y)} \cdot \frac{V(X)}{V(X^*)} \cdot \frac{V(Y)}{V(Y^*)}\\
&-\left( \frac{V(Y)-V(Y^*)}{V(Y^*)} + \frac{V(X)-V(X^*)}{V(X^*)} + \frac{V(X)-V(X^*)}{V(X^*)}\cdot \frac{V(Y)-V(Y^*)}{V(Y^*)} \right)\\
&= d(X,Y) \cdot \frac{V(X)}{V(X^*)} \cdot \frac{V(Y)}{V(Y^*)} \\
&- \left( \frac{V(Y)}{V(Y^*)} - 1 \right) - \left( \frac{V(X)}{V(X^*)} - 1 \right) - \left( \frac{V(X)}{V(X^*)} - 1 \right) \left( \frac{V(Y)}{V(Y^*)} - 1 \right)\\
&= d(X,Y) \cdot \frac{V(X)}{V(X^*)} \cdot \frac{V(Y)}{V(Y^*)} + 1 - \frac{V(X)}{V(X^*)} \cdot \frac{V(Y)}{V(Y^*)}\\
&= \frac{V(X)}{V(X^*)} \cdot \frac{V(Y)}{V(Y^*)} \cdot \left[ d(X,Y) - 1 \right] + 1.
\end{align*}
Since $d(X,Y) \leq 1$, we have $d(X,Y) - 1 \leq 0$. Also, $\frac{V(X)}{V(X^*)} \le \frac{1}{1-\delta}$ and $\frac{V(Y)}{V(Y^*)} \le \frac{1}{1-\delta}$, so
\[
d(X^*,Y^*) \geq \left( \frac{1}{1-\delta} \right)^2 \left[ d(X,Y) - 1 \right] + 1,
\]
which gives a lower bound for $d(X^*,Y^*)$. Together with the upper bound (1), we have
\[
\left( \frac{1}{1-\delta} \right)^2 \left[ d(X,Y) - 1 \right] + 1 \leq d(X^*,Y^*) \leq \left(\frac{1}{1-\delta}\right)^2 \cdot d(X,Y).
\]
Subtracting $d(X,Y)$ throughout,
\[
\left( \left( \frac{1}{1-\delta} \right)^2 - 1 \right) \left[ d(X,Y) - 1 \right] \leq d(X^*,Y^*) - d(X,Y) \leq \left( \left( \frac{1}{1-\delta} \right)^2 - 1 \right) d(X,Y).
\]
Observe that the left side is nonpositive and the right side is nonnegative, so
\[
| d(X^*,Y^*) - d(X,Y) | \leq \left( \left( \frac{1}{1-\delta} \right)^2 - 1 \right) \cdot \max \{ d(X,Y), 1 - d(X,Y) \}.
\]
It's easy to check that $\left( \left( \frac{1}{1-\delta} \right)^2 - 1 \right) < 4 \delta$ for $0 < \delta < \frac{1}{4}$, and $0 \leq d(X,Y) \leq 1$, so
\[
| d(X^*,Y^*) - d(X,Y) | < 4 \delta \cdot 1 = 4 \delta,
\]
which proves the theorem. $\square$

\begin{theorem}[continuity of the discrepancy]
Let $\A$ be an $m \times n$ array with row set $R$ and column set $C$, now
\[
\sum_{i \in R} \sum_{j \in C} a_{ij} = 1, \qquad a_{ij} \ge 0 \quad \forall i, j
\]
We assume that there is no dominant row or dominant column; precisely, assume that $\frac{c_1}{m} \le d_{\text{row},i} \le \frac{c_2}{m}$ and $\frac{c_3}{n} \le d_{\text{col},i} \le \frac{c_4}{n}$ for some constants $c_1,c_3 \leq 1$ and $c_2,c_4 \geq 1$. Suppose $X^* \subset X \subset R$ and $Y^* \subset Y \subset C$, and also suppose $\frac{\Vol(X^*)}{\Vol(X)} \ge 1-\delta$ and $\frac{\Vol(Y^*)}{\Vol(Y)} \ge 1-\delta$ for some $\delta < \frac{1}{4}$. 
Then
\[
\disc(X^*,Y^*;X,Y) \le 4\delta\sqrt{\frac{c_2 c_4}{c_1 c_3}}\,.
\]
\end{theorem}

\noindent
\textbf{Proof.}
By the definitions of $\Vol(X)$ and $\Vol(Y)$,
\[
\Vol(X) = \sum_{i \in X} d_{\text{row},i} \ge \frac{c_1}{m}\cdot|X|,
\]
so $\frac{1}{\Vol(X)} \le \frac{m}{c_1}\cdot\frac{1}{|X|}$, and similarly
\[
\Vol(Y) = \sum_{j \in Y} d_{\text{col},j} \ge \frac{c_3}{n}\cdot|Y|,
\]
so $\frac{1}{\Vol(Y)} \le \frac{n}{c_3}\cdot\frac{1}{|Y|}$. Note that $|X|$ is the number of rows in $X$, and $|Y|$ is the number of columns in $Y$.

\vspace{0.1in}
\noindent Since $\frac{c(X,Y)}{\Vol(X)} \le 1$ and $\frac{c(X,Y)}{\Vol(Y)}\le 1$,
\[
\rho(X,Y)= \frac{c(X,Y)}{\Vol(X)\Vol(Y)} \le \frac{1}{\Vol(X)} \le \frac{m}{c_1\cdot |X|} \,.
\] 
Likewise,
\[
\rho(X,Y) \le \frac{1}{\Vol(Y)} \le \frac{n}{c_3\cdot |Y|} \,.
\]
Hence,
\[
\rho(X,Y) \le \min \left\{ \frac{m}{c_1\cdot |X|}\,, ~\frac{n}{c_3\cdot |Y|} \right\} .
\]
Let $K = \min \left\{ \frac{m}{c_1\cdot |X|},\frac{n}{c_3\cdot |Y|} \right\}$.


\vspace{0.1in}
\noindent To find the upper bound for $\rho(X^*,Y^*)-\rho(X,Y)$.
\begin{align*}
\rho(X^*,Y^*) &= \frac{c(X^*,Y^*)}{\Vol(X^*)\Vol(Y^*)}\\
&\le \frac{c(X,Y)}{\Vol(X^*)\Vol(Y^*)}\\
&= \frac{c(X,Y)}{\Vol(X)\Vol(Y)}\cdot \frac{\Vol(X)\Vol(Y)}{\Vol(X^*)\Vol(Y^*
)}\\
&\le \rho(X,Y) \cdot \frac{1}{(1-\delta)^2}\, .
\end{align*}
This gives the upper bound as
\begin{equation} \label{upper-bound}
\rho(X^*,Y^*)-\rho(X,Y) \le \rho(X,Y)\cdot \left(\frac{1}{(1-\delta)^2}-1 \right).
\end{equation}
Next, to find the lower bound, we know that 
\[
c(X^*,Y^*) = c(X,Y)-[c(X-X^*,Y^*)+c(X^*,Y-Y^*)+c(X-X^*,Y-Y^*)].
\]
Since
\begin{align*}
c(X-X^*,Y^*) &\le K\cdot (\Vol(X)-\Vol(X^*)) \cdot \Vol(Y^*),\\
c(X^*,Y-Y^*) &\le K\cdot \Vol(X^*) \cdot (\Vol(Y)-\Vol(Y^*)),\\
c(X-X^*,Y-Y^*) &\le K\cdot (\Vol(X)-\Vol(X^*)) \cdot (\Vol(Y)-\Vol(Y^*)),
\end{align*}
it follows that
\[
c(X^*,Y^*) \ge c(X,Y)-K\cdot \Vol(X)\Vol(Y) + K \cdot \Vol(X^*)\Vol(Y^*) .
\]

\vspace{0.1in}
\noindent Dividing both sides by $\Vol(X^*)\Vol(Y^*)$,
\begin{align*}
\rho(X^*,Y^*) &\ge \frac{c(X,Y)}{\Vol(X^*)\Vol(Y^*)}-K\cdot \frac{\Vol(X)\Vol(Y)}{\Vol(X^*)\Vol(Y^*)} + K\\
&= \rho(X,Y) \cdot \frac{\Vol(X)\Vol(Y)}{\Vol(X^*)\Vol(Y^*)} -K\cdot \frac{\Vol(X)\Vol(Y)}{\Vol(X^*)\Vol(Y^*)} + K\\
&= (\rho(X,Y) - K) \cdot \frac{\Vol(X)\Vol(Y)}{\Vol(X^*)\Vol(Y^*)} + K.
\end{align*}

\vspace{0.1in}
\noindent Subtracting $\rho(X,Y)$ on both sides gives the lower bound as
\begin{align}
\rho(X^*,Y^*)-\rho(X,Y) &\ge (\rho(X,Y) - K) \cdot \frac{\Vol(X)\Vol(Y)}{\Vol(X^*)\Vol(Y^*)} - [\rho(X,Y) - K]\\
&= (\rho(X,Y)-K) \left( \frac{\Vol(X)\Vol(Y)}{\Vol(X^*)\Vol(Y^*)}-1 \right)\\
&= (\rho(X,Y)-K) \left( \frac{1}{(1-\delta)^2}-1 \right).  \label{lower-bound}
\end{align}
Together, the upper bound \eqref{upper-bound} and the lower bound \eqref{lower-bound} give
\[
(\rho(X,Y)-K) \cdot \left( \frac{1}{(1-\delta)^2}-1 \right) \le \rho(X^*,Y^*)-\rho(X,Y) \le \rho(X,Y)\cdot \left( \frac{1}{(1-\delta)^2}-1 \right),
\]
which implies
\[
|\rho(X^*,Y^*)-\rho(X,Y)| \le \left( \frac{1}{(1-\delta)^2}-1 \right) \cdot \max \{K-\rho(X,Y), ~\rho(X,Y)\}.
\]
We have $\frac{1}{(1-\delta)^2}-1 \le 4\delta$ (since $\delta < \frac{1}{4}$) and $\max\{K-\rho(X,Y),\rho(X,Y)\} \le K$, so
\[
|\rho(X^*,Y^*)-\rho(X,Y)| \le 4 \delta \cdot K
\]
which gives the bound of
\[
\disc(X^*,Y^*;X,Y) \le 4\delta \cdot K \cdot \sqrt{\Vol(X^*)\Vol(Y^*)}
\]
on the discrepancy. Using the fact that $\min\{a,b\}\le\sqrt{ab}$ for any $a,b \ge 0$,
\begin{equation} \label{K-bound}
K = \min \left\{\frac{m}{c_1 |X|},~\frac{n}{c_3 |Y|} \right\} \le \sqrt{\frac{m}{c_1 |X|} \cdot \frac{n}{c_3 |Y|}}\,.
\end{equation}
By the definitions of $\Vol(X^*)$ and $\Vol(Y^*)$,
\begin{align}
\Vol(X^*) &= \sum_{i \in X^*} d_{\text{row},i} \le \frac{c_2}{m}\cdot|X^{*}|, \label{Vol-X*-bound}\\
\Vol(Y^*) &= \sum_{j \in Y^*} d_{\text{col,}j} \le \frac{c_4}{n}\cdot|Y^{*}|. \label{Vol-Y*-bound}
\end{align}
By \eqref{K-bound}, \eqref{Vol-X*-bound}, and \eqref{Vol-Y*-bound}, it follows that
\begin{align*}
\disc(X^*,Y^*;X,Y) &\le 4\delta \cdot \sqrt{\frac{m}{c_1\cdot |X|} \cdot \frac{n}{c_3\cdot |Y|}}\cdot \sqrt{\frac{c_2}{m}\cdot|X^{*}|\cdot \frac{c_4}{n}\cdot|Y^{*}|}\\
&= 4\delta \sqrt{\frac{c_2 c_4}{c_1 c_3} \cdot \frac{|X^*|}{|X|} \cdot \frac{|Y^*|}{|Y|}}.
\end{align*}
Since $\frac{|X^*|}{|X|}\le 1$ and $\frac{|Y^*|}{|Y|}\le 1$,
\[
\disc(X^*,Y^*;X,Y) \le 4\delta \sqrt{\frac{c_2 c_4}{c_1 c_3}}\,,
\]
which proves the theorem. $\square$

\section{Some other properties of the discrepancy}
\label{prop}

\begin{proposition} Assume $\A$ is an independent table, that is, $a_{ij} = 
d_{row,i} d_{col,j}$ for all $i,j$. Assume that 
there are not identically zero rows or columns (and hence, $a_{ij}>0$
for all $i,j$). Then $\disc_k(A) = 0$ for $1 \leq k \leq \rk(A)$.
\end{proposition}

\noindent
\textbf{Proof.}
It is trivial, because $\rho (X,Y) =1$  for all $X\subset R$, $Y\subset C$.
$\square$

\begin{proposition} Let $\A$ be a table with row partition $(R_1, \ldots, R_k)$ and column partition $(C_1, \ldots, C_k)$. Assume that $a_{ij} \geq 0$. Define $s_{ab} = \sum_{i \in R_a} \sum_{j \in C_b} a_{ij}$ for all $1 \leq a,b \leq k$, and suppose that  the $k\times k$ matrix $(s_{ab})$ defines an independent table. 
Then $\disc_k(\A ), \dots ,\disc_2 (\A ) \leq \disc_1(\A )$. 
We call such a table $k$-way contracted independent.
\end{proposition}

\noindent
\textbf{Proof.}
It is trivial, because with the above partitions, $\rho (R_a , C_b) =1$, so
$\disc_k (X, Y; R_a ,C_b ) =\disc_k (X, Y; R ,C )$ for all 
$X\subset R_a$, $Y\subset C_b$ $(a,b=1,\dots ,k)$.

Note that if $A$ is $k$-way contracted independent with some $k$, than it is 
also $\ell$-way contracted independent with $\ell =2,3,\dots ,k-1$.
$\square$

\begin{proposition} Let $A$ be a table with row partition $(R_1, \ldots, R_k)$ and column partition $(C_1, \ldots, C_k)$. Assume  that for all $1 \leq a,b \leq k$, $i \in R_a$, and $j \in C_b$, we have $a_{ij} = c_{ab}$, where $c_{ab}>0$ depends only on $a,b$. Then $0 = \disc_k(\A) = \disc_{k+1}(\A) = \cdots = 
\disc_{\rk \A}(\A)$.
\end{proposition}

\noindent
\textbf{Proof.}
It is trivial, because with the above partitions, 
$\rho (X, Y) =\rho ( R_a ,C_b )$ for all 
$X\subset R_a$, $Y\subset C_b$ $(a,b=1,\dots ,k)$.
$\square$

\begin{lemma}
For any integer $n>1$ and arbitrary positive real numbers $u_1 ,\dots ,u_n$ and
$v_1 ,\dots v_n$ we have
$$
 \min_{1\le i \le n} \frac{u_i}{v_i} \le
 \frac{u_1 +\dots +u_n}{v_1 +\dots +v_n }\le\max_{1\le i \le n} \frac{u_i}{v_i},
$$
and equality holds if and only if the ratios $\frac{u_i}{v_i}$ have the
same value.
\end{lemma}

\begin{proposition} Let $A$ be a table of nonnegative entries,
$X\subset R$, $Y\subset C$. With $2\le k \le |X|$ let $X_1 ,\dots ,X_k$
be a proper $k$-partition of $X$. Then
$$
 \min_i \rho (X_i ,C) \le \rho (X ,C) \le  \max_i \rho (X_i ,C) .
$$
The same holds for $k$-partitions of $Y$ too.
\end{proposition}

\noindent
\textbf{Proof.} It is straightforward with the Lemma.
$\square$

A bit more is true for 0-1 tables: the density of the union of two disjoint
stripes is the weighted average of the parts,  in the proportion of their
volumes.
\begin{proposition}
Let $A$ be the adjacency matrix of a graph, $X,X' \subset V$ are disjoint,
$Y\subset V$ is arbitrary. Assume that $V (X' ) =v V(X)$, $v>0$. Then
$$
 d(X \cup X' ,Y) = \frac{1}{1+v} d(X,Y )+ \frac{v}{1+v} d(X',Y ) .
$$
\end{proposition}

\noindent
\textbf{Proof.}
The proof is trivial using that $e(X \cup X' ,Y)= e(X,Y) +e(X',Y)$ and the
volumes are also added together. If $v=1$, we get the arithmetic average.
$\square$

Some other remarks.
\begin{itemize}
\item
Note that $\A$ is $2\times 2$ contracted independent if and only if 
(after possibly permuting its rows and columns) there is a 2-partition
of the rows and that of the columns such that the contracted matrix
$$
 \begin{pmatrix} s_{11} & s_{12} \\
                 s_{21} & s_{22}
 \end{pmatrix}
$$
satisfies $s_{11} s_{22} =s_{12} s_{21}$.

With some continuity arguments (there are no dominant rows and columns) we
can get an approximate solution. In the graph case ($\A$ is symmetric),
the minimum and maximum cuts are the two extremes.
\item
If the table is $k\times k$ contracted independent, the above equation has
$k-1$ different solutions. We conjecture that partitions producing 
such tables are good candidates for minimizing the $k$-way discrepancy.

\item
We may consider the $k$-partition of the rows and columns such that
$$
 \disc_k (\A ) = 
 \disc (\A ; R_1 ,\dots ,R_k , C_1 ,\dots ,C_k ).
$$
Let $a,b$ and $X\subset R_a$, $Y\subset C_b$ be such that 
$$
 \disc (X,Y; R_a ,C_b )=
\disc (\A ; R_1 ,\dots ,R_k , C_1 ,\dots ,C_k ).
$$
Then try to divide $R_a$ into $X, {\bar X}$ and $C_b$ into $Y, {\bar Y}$ and
consider the so obtained $(k+1)$-partitions.
\end{itemize}
 


\section*{Acknowledgements}

The research was supported by 
the Budapest Semester of Mathematics, undergraduate research program 
(RES5, 2015).


\begin{thebibliography}{1}

\bibitem{Bolla14} 
Bolla, M., SVD, discrepancy, and regular structure of contingency tables,
\textit{Discrete Applied Mathematics} \textbf{176} (2014), 3-11.


\bibitem{Bollarx}
Bolla, M., Relating multiway discrepancy and singular values of 
nonnegative rectangular matrices, 
\textit{Discrete Applied Mathematics} \textbf{203} (2016), 26-34.

\end{thebibliography}
\end{document}